\documentclass[10pt, a4paper]{amsart}
\usepackage{amsfonts}
\usepackage{amssymb,amsthm, amsmath}
\usepackage{bbm}
\usepackage{url}

\newtheorem{theorem}{Theorem}
\newtheorem{definition}[theorem]{Definition}
\newtheorem{proposition}[theorem]{Proposition}
\newtheorem{corollary}[theorem]{Corollary}
\newtheorem{lemma}[theorem]{Lemma}

\theoremstyle{remark}

\title{Indispensable binomials in semigroup ideals}
\author{Ignacio Ojeda \and
Alberto Vigneron-Tenorio}

\address{Departamento de Matem\'aticas, Universidad de Extremadura,
E-06071 Badajoz (Spain).} \email{ojedamc@unex.es}

\address{Departamento de Matem\'aticas, Universidad de C\'adiz,
E-11405 Jerez de la Frontera (Spain).} \email{alberto.vigneron@uca.es}



\begin{document}

\date{\today}

\begin{abstract}
In this paper, we deal with the problem of uniqueness of minimal system of binomial generators of a semigroup ideal. Concretely, we give different necessary and/or sufficient conditions for uniqueness of such minimal system of generators. These conditions come from the study and combinatorial description of the so-called indispensable binomials in the semigroup ideal.
\end{abstract}

\maketitle

\section*{Introduction}

Under suitable conditions on a finitely generated semigroup $\mathcal{A} = \mathbf{a}_1 \mathbb{N} + \ldots + \mathbf{a}_r \mathbb{N},$ a binomial ideal $I_\mathcal{A}$ of the polynomial ring in $r$ variables over a field $\mathbbmss{k},\ R = \mathbbmss{k}[X_1, \ldots, X_r],$ is determined: the so-called semigroup ideal of $\mathcal{A}$
(see Section \ref{Sect 01} for the details).

Semigroup ideals play a relevant role in (Computational) Commutative Algebra and Algebraic Geometry. Moreover, they have a lot of applications in different areas such as Statistics, Operational Research or Phylogenetic reconstruction (to cite only three examples).

As in the case of monomial ideals, semigroup ideals have a rich combinatorial structure coming directly from the semigroup (see, e.g. \cite{Collectanea, Charalambous08} and \cite{OjVi}). However, although reasonable conditions on $\mathcal{A}$ guarantee the existence of minimal systems of binomial generators for $I_\mathcal{A},$ they need not be unique, in contrast with the monomial case.

Many successful efforts have been made to compute and describe (minimal) systems of binomial generators of semigroup ideals (see, e.g., \cite{BCMP} or Chapter 9 in \cite{Rosales} and the references therein), but only recently the problem of the uniqueness has been explicitly treated. In fact, the interest in this problem arises first from Algebraic Statistics (\cite{Takemura}) in which the systems of binomial generators of particular families of semigroup ideals defining certain statistics models are called Markov bases of the model (\cite{Diaconis}). Subsequently, different authors have investigated this uniqueness problem. In \cite{Aoki} and in \cite{Ohsugi}, the notions of indispensable binomial and monomial are introduced, respectively, and in \cite{Charalambous07} both properties are studied from a combinatorial point of view. It is convenient to observe that one can find in the literature results which imply a characterization of particular families of semigroups ideals with a unique system of binomial generators (see Corollary 3.4 in \cite{Peeva}, Corollary 3.11 in \cite{OjPis} or Proposition 1 in \cite{PisVi}).
Furthermore, it is known that the generic lattice ideals in introduced by I. Peeva and B. Sturmfels in \cite{Peeva}
have a unique minimal system of binomial generators (see Lemma 3.3 and Remark 4.4.3 in \cite{Peeva}); therefore, one may assure that there exist semigroup ideals with arbitrary large systems of indispensable binomials generators (see \cite{Oj2}). All these facts suggest that the uniqueness problem might be connected with other interesting problems about semigroup ideals.

Recall that a binomial in a semigroup ideal  $I_\mathcal{A}$ is said to be indispensable if it belongs (up to a scalar multiple) to every system of binomial generators $I_\mathcal{A}$ (see \cite{Ohsugi05}). Thus, $I_\mathcal{A}$ has a unique system of binomial generators if, and only if, it is generated by indispensable binomials. On the other hand, a monomial is called indispensable if it appears (up to a scalar multiple) in, at least, one binomial of any system of binomial generators of $I_\mathcal{A}.$

Indispensable monomials were introduced by Aoki, Takemura and
Yoshida in \cite{Aoki}. They are considered as a first approximation to the study of indispensable binomials, notice that any indispensable binomial is a difference of two indispensable monomials. Moreover, indispensable monomials always exist (see Proposition 3.1 in \cite{Charalambous07}) in clear contrast to the indispensable binomials (consider, e.g. $I_\mathcal{A} = \langle x-y,x-z \rangle$).

In the first part of this paper (Section \ref{Sect 02}), we give a combinatorial necessary and sufficient condition for the existence of indispensable binomials in a semigroup ideal (Theorem \ref{Th Indis}). Our condition depends on the knowledge of some simplicial complexes associated to the semigroup introduced by S. Eliahou in his PhD Thesis (\cite{Eliahou}). As a consequence, an explicit characterization of all indispensable binomials and monomials of a semigroup ideal $I_\mathcal{A}$ is given. In Section \ref{Sect 03}, we deal with the problem of the existence of indispensable binomials in a semigroup ideal $I_\mathcal{A}$ by using only Gr\"obner bases techniques; concretely, we give an effective necessary and sufficient condition for the existence of indispensable binomials in $I_\mathcal{A}$ consisting in the computation of $r$ (the number of indeterminates of the corresponding polynomial ring) Gr\"obner bases (Theorem \ref{Th indisGB}). We end the paper by illustrating our results with an example borrowed from  Algebraic Statistics.

\section{Preliminaries: Definitions and notation}\label{Sect 01}

Let $\mathcal{A}$ be a commutative semigroup with zero element $0 \in \mathcal{A}$ and let $G(\mathcal{A})$ be a commutative group with a semigroup homomorphism $\iota : \mathcal{A} \rightarrow G(\mathcal{A})$ such that every homomorphism from $\mathcal{A}$ to a group factors in a unique way through $\iota.$ The commutative group $G(\mathcal{A})$ exists and is unique up to isomorphism; it is called the associated commutative group of $\mathcal{A}.$ Further, $G(\mathcal{A})$ is finitely generated when $\mathcal{A}$ is. The map $\iota$ is injective if, and only if, $\mathcal{A}$ is cancellative, that is to say, if $\mathbf{a} + \mathbf{b} = \mathbf{a} + \mathbf{c},\ \mathbf{a}, \mathbf{b}, \mathbf{c} \in \mathcal{A},$ implies $\mathbf{b} = \mathbf{c};$ in this case, $G(\mathcal{A})$ is the smallest group containing $\mathcal{A}.$

Through all this paper, $\mathcal{A}$ will denote a finitely generated commutative cancellative semigroup with zero element. Moreover, we will assume always that $\mathcal{A}$ is combinatorially finite, that is to say, there are only finitely many ways to write $\mathbf{a} \in \mathcal{A} \setminus \{0\}$ as a sum $\mathbf{a} = \mathbf{a}_1 + \ldots + \mathbf{a}_q,$ with $\mathbf{a}_i \in \mathcal{A} \setminus \{0\}.$ Equivalently, $\mathcal{A}$ is combinatorially finite if, and only if, $\mathcal{A} \cap (-\mathcal{A}) = \{0\}$ (see Proposition 1.1 in  \cite{BCMP}). Notice that this property guarantees that $\mathbf{b} \preceq_\mathcal{A} \mathbf{a} \Longleftrightarrow \mathbf{a}-\mathbf{b} \in \mathcal{A}$ is a well defined partial order on $\mathcal{A}.$

Examples of finitely generated commutative cancellative combinatorially finite semigroups with zero element (semigroups in the following) are the subsemigroup of $\mathbb{N}^d$ generated by the columns of a non-negative integer matrix $A \in \mathbb{N}^{d \times r}.$

\medskip
Let $\mathbbmss{k}$ be a fixed arbitrary field. We write $\mathbbmss{k}[\mathcal{A}]$ for the $\mathbbmss{k}$-vector space $$\mathbbmss{k}[\mathcal{A}] = \bigoplus_{\mathbf{a} \in \mathcal{A}} \mathbbmss{k}\, \mathbf{t}^\mathbf{a}$$ endowed with a multiplication which is $\mathbbmss{k}$-linear and such that $\mathbf{t}^\mathbf{a} \cdot \mathbf{t}^\mathbf{b} := \mathbf{t}^{\mathbf{a}+\mathbf{b}},\ \mathbf{a}, \mathbf{b} \in \mathcal{A}.$ Thus $\mathbbmss{k}[\mathcal{A}]$ has a natural $\mathbbmss{k}$-algebra structure and we will refer to it as the semigroup algebra of $\mathcal{A}.$

The choice of a system of generators $\mathbf{a}_1, \ldots, \mathbf{a}_r$ of $\mathcal{A}$
induces a natural $\mathcal{A}-$grading on $R := \mathbbmss{k}[X_1, \ldots, X_r],$ by assigning weight $\mathbf{a}_i$ to $X_i,\ i = 1, \ldots, r;$ that is to say, $$R = \bigoplus_{\mathbf{a} \in \mathcal{A}} R_\mathbf{a},$$ where $R_\mathbf{a}$ is the vector subspace of $R$ generated by all the monomials $X^\mathbf{u} := X_1^{u_1} \cdots X_r^{u_r}$ with $\sum_{i=1}^r u_i \mathbf{a}_i = \mathbf{a}$ and $\mathbf{u} = (u_1, \ldots, u_r)' \in \mathbb{N}^r$ (the dash means transpose). Since $\mathcal{A}$ is combinatorially finite, the vector spaces $R_\mathbf{a}$ are finite dimensional (see Proposition 1.2 in \cite{BCMP}). Concretely, there are $\dim_\mathbbmss{k} R_\mathbf{a}$ monomials of $\mathcal{A}$-degree $\mathbf{a},$ for each $\mathbf{a} \in \mathcal{A}.$ Let $M_\mathbf{a}$ denote the set of monomials in $R_\mathbf{a}.$

In the following, we will assume that $\mathbf{a}_1, \ldots, \mathbf{a}_r$ is a fixed system of generators of $\mathcal{A}.$

The natural semigroup morphism $\pi : \mathbb{N}^r \to \mathcal{A};\ \mathbf{u} \mapsto \sum_{i=1}^r u_i \mathbf{a}_i$ defines a $\mathcal{A}$-(multi)graded surjective $\mathbbmss{k}$-algebra morphism $$\varphi_0 : R = \mathbbmss{k}[\mathbb{N}^r] \longrightarrow \mathbbmss{k}[\mathcal{A}];\ X_i \longmapsto \mathbf{t}^{\mathbf{a}_i}.$$ Thus, the ideal $I_\mathcal{A} := \ker(\varphi_0)$ is a  $\mathcal{A}$-homogeneous ideal called the (semigroup) ideal of $\mathcal{A}.$ Notice that $I_\mathcal{A}$ is a toric ideal when $G(\mathcal{A})$ is torsion free (see \cite{Sturmfels95}, chapter 4).

It is well known (see \cite{Herzog70}) that $I_\mathcal{A}$ is the ideal of $R$ generated by $$\Big\{X^\mathbf{u} - X^\mathbf{v} : \sum_{i=1}^r u_i \mathbf{a}_i = \sum_{i=1}^r v_i \mathbf{a}_i \Big \},$$ where $X^\mathbf{u} = X_1^{u_1} \cdots X_r^{u_r}$ and $X^\mathbf{v} = X_1^{v_1} \cdots X_r^{v_r},$ as usual. Therefore, there exist minimal systems of $\mathcal{A}$-homogeneous generators of $I_\mathcal{A}$ consisting in finitely many pure difference binomials, i.e. differences of two monomials of the same $\mathcal{A}-$degree (see, e.g. \cite{BCMP}, Section 2). In fact, one has the following:

\begin{lemma}\label{Lemma1}
Every system of binomial generators of $I_\mathcal{A}$ is $\mathcal{A}-$graded.
\end{lemma}

\begin{proof}
By definition $I_\mathcal{A} = \ker(\varphi_0).$ So, if $X^\mathbf{u} - \lambda X^\mathbf{v} \in I_\mathcal{A},$ then $$\mathbf{t}^{\sum_{i=1}^r u_i \mathbf{a}_i} = \lambda \mathbf{t}^{\sum_{i=1}^r v_i \mathbf{a}_i},$$ that is to say, $\sum_{i=1}^r u_i \mathbf{a}_i = \sum_{i=1}^r v_i \mathbf{a}_i$ and $\lambda = 1.$
\end{proof}

The $\mathcal{A}-$degrees of the polynomials appearing in any minimal system of $\mathcal{A}$-homogeneous generators of $I_\mathcal{A}$ do not depend on the system of generators: it is well known that \emph{the number of polynomials of $\mathcal{A}-$degree $\mathbf{a} \in \mathcal{A}$ in a minimal system of $\mathcal{A}$-homogeneous generators is $\dim_\mathbbmss{k} \mathrm{Tor}_1^R(\mathbbmss{k}, \mathbbmss{k}[\mathcal{A}])_\mathbf{a}$} (see, e.g. \cite{Miller05}, Section 8.3). Thus, we say that $I_\mathcal{A}$ has minimal generators in degree $\mathbf{a}$ when $\dim_\mathbbmss{k} \mathrm{Tor}_1^R(\mathbbmss{k}, \mathbbmss{k}[\mathcal{A}])_\mathbf{a} \neq 0.$

\begin{definition}
We say that $\mathbf{a} \in \mathcal{A}$ is a minimal $\mathcal{A}-$degree of $I_\mathcal{A},$ if $I_\mathcal{A}$ has minimal generators in degree $\mathbf{a}.$
\end{definition}

Note that there are finitely many minimal $\mathcal{A}-$degrees of $I_\mathcal{A}.$ 

\begin{definition}\label{Def Indis}
We say that $\mathbf{a} \in \mathcal{A}$ is an indispensable $\mathcal{A}-$degree of $I_\mathcal{A},$ if
every (minimal) system of $\mathcal{A}$-homogeneous generators of $I_\mathcal{A}$ contains one, and only one, polynomial of $\mathcal{A}-$degree $\mathbf{a}.$  In this case, such polynomial is a binomial and it is said to be an indispensable binomial of $I_\mathcal{A}.$
\end{definition}

Notice first that every indispensable $\mathcal{A}-$degree is a minimal $\mathcal{A}-$degree. Furthermore, observe that an indispensable binomial of $I_\mathcal{A}$ appears (up to a scalar multiple) in every system of $\mathcal{A}$-homogeneous generators of $I_\mathcal{A}.$ Moreover, it is easy to see that a binomial that appears (up to scalar multiple) in every system of $\mathcal{A}$-homogeneous generators of $I_\mathcal{A}$ is indispensable. Indeed, if $\mathcal{S}$ is a system of $\mathcal{A}$-homogeneous generators of $I_\mathcal{A}$ containing two different polynomials, $f$ and $g,$ of the same $\mathcal{A}-$degree, by substituting $f$ by $f-g$ in $\mathcal{S},$ we obtain a new system of $\mathcal{A}$-homogeneous generators of $I_\mathcal{A}$ not containing $f.$ That is to say, if $f \in I_\mathcal{A}$ is not indispensable, it does not belong to every system of $\mathcal{A}$-homogeneous generators of $I_\mathcal{A}.$

Summarizing, we have that:

\begin{proposition}
The following statements are equivalent:
\begin{itemize}
\item[(a)] $I_\mathcal{A}$ has a unique (up to a scalar multiple) minimal system of $\mathcal{A}$-homogeneous generators.
\item[(b)] Every minimal system of $\mathcal{A}$-homogeneous generators consists of binomials.
\item[(c)] $I_\mathcal{A}$ has a unique minimal system of (pure difference) binomial generators.
\item[(d)] $I_\mathcal{A}$ is generated by its indispensable binomials.
\end{itemize}
\end{proposition}

\section{Combinatorial description of indispensability}\label{Sect 02}

In this section, we will give a necessary and sufficient condition for the existence of indispensable binomials in $I_\mathcal{A}$ for a given semigroup $\mathcal{A} = \mathbf{a}_1 \mathbb{N} + \ldots + \mathbf{a}_r \mathbb{N}.$

\medskip
We start by introducing a combinatorial object associated to $\mathcal{A}.$

\begin{definition}
For any $\mathbf{a} \in \mathcal{A},$ define the abstract simplicial complex $\nabla_\mathbf{a}$ on the vertex set $M_\mathbf{a} = \{X^\mathbf{u} = X_1^{u_1} \cdots X^{u_r} \mid \sum_{i=1}^r u_i \mathbf{a}_i = \mathbf{a}\}$ $$\nabla_\mathbf{a} = \{ F \subseteq M_\mathbf{a} \mid \gcd(F) \neq 1 \},$$ where $\gcd(F)$ denotes the greatest common divisor of the monomials in $F.$
\end{definition}

Recall that $M_\mathbf{a}$ has finitely many monomials because of combinatorially finiteness of $\mathcal{A}.$ So, the simplicial complexes $\nabla_\mathbf{a}$ are finite.

\medskip
The next proposition was first proved by S. Eliahou who introduced the simplicial
complexes $\nabla_\mathbf{a}$ in \cite{Eliahou}. Other proofs, in a more
general context, can be found in \cite{OjVi} and in \cite{Charalambous08}.

\begin{proposition}\label{Prop EliahouGeneralized}
Let $\mathbf{a} \in \mathcal{A}.\ \mathbf{a}$ is a minimal $\mathcal{A}-$degree of $I_\mathcal{A}$ if,
and only if, $\nabla_\mathbf{a}$ is not connected.
\end{proposition}

It is important to observe that the $1-$skeleton of $\nabla_\mathbf{a}$ is a subgraph of the
graph given in \cite{Charalambous07} (Definition 2.1) with the same set of vertices and the same number connected components. Thus, one can use the simplicial complexes $\nabla_\mathbf{a}$ to obtain the same results as in that paper. In fact, Theorems 2.6 and 2.7 in \cite{Charalambous07} can be understood as a new version of Theorem 2.5 in \cite{BCMP} for the simplicial complexes $\nabla_\mathbf{a},$ by taking into account that both the simplicial complexes in \cite{BCMP} and the simplicial complexes $\nabla_\mathbf{a}$ have isomorphic homology $\mathbbmss{k}-$vector spaces (see \cite{OjVi}, Theorem 3).

\medskip
An immediate consequence of Proposition \ref{Prop EliahouGeneralized} is the following:

\begin{corollary}\label{Cor Indis}
Let $\mathbf{a} \in \mathcal{A}.\ \mathbf{a}$ is an indispensable $\mathcal{A}-$degree of $I_\mathcal{A}$ if, and only if, $\nabla_\mathbf{a} = \Big\{ \{X^\mathbf{u}\}, \{X^\mathbf{v}\} \Big\}.$
\end{corollary}

\begin{proof}
Let $\mathbf{a}$ be an indispensable $\mathcal{A}-$degree of $I_\mathcal{A}$ and let $X^\mathbf{u} - X^\mathbf{v} \in I_\mathcal{A}$ be the corresponding indispensable binomial. If there exists a monomial $X^\mathbf{w} \in M_\mathbf{a}$ different from $X^\mathbf{u}$ and $X^\mathbf{v},$ we could replace $X^\mathbf{u} - X^\mathbf{v}$ by $X^\mathbf{u} - X^\mathbf{w}$ and $X^\mathbf{v} - X^\mathbf{w}$ in a (minimal) system of generators of $I_\mathcal{A},$ thus obtaining a (not-necessarily minimal) system of generators of $I_\mathcal{A}$ not containing $X^\mathbf{u} - X^\mathbf{v}$ which is not possible by definition \ref{Def Indis}. Therefore, $M_\mathbf{a} = \{X^\mathbf{u}, X^\mathbf{v}\}$ and, by Proposition \ref{Prop EliahouGeneralized}, we conclude that $\nabla_\mathbf{a} = \Big\{ \{X^\mathbf{u}\}, \{X^\mathbf{v}\} \Big\},$ i.e., $\gcd\big(X^\mathbf{u}, X^\mathbf{v}\big) = 1$ (recall that every indispensable $\mathcal{A}-$degree is a minimal $\mathcal{A}-$degree).

Conversely, if $\nabla_\mathbf{a} = \Big\{ \{X^\mathbf{u}\}, \{X^\mathbf{v}\} \Big\},$ by Proposition \ref{Prop EliahouGeneralized}, $\mathbf{a}$ is a minimal $\mathcal{A}-$degree of $I_\mathcal{A}.$ Moreover, since the only polynomial of $\mathcal{A}-$degree (up to scalar multiple) is $X^\mathbf{u} - X^\mathbf{v},$ we conclude that it has to be indispensable and so $\mathbf{a}$ is an indispensable $\mathcal{A}-$degree of $I_\mathcal{A}.$
\end{proof}

The above result was also noticed by H. Charalambous et al. (see Theorem 4.1 in \cite{Charalambous07}).

%
%
%

\begin{theorem}\label{Th Indis}
Let $\mathbf{a} \in \mathcal{A}.\ \mathbf{b} \prec_\mathcal{A} \mathbf{a}$ is indispensable if, and only if, there exists $\{X^\mathbf{u}, X^\mathbf{v}\} \in \nabla_\mathbf{a}$ such that
\begin{itemize}
\item[(a)] $\gcd(X^\mathbf{u}, X^\mathbf{v}) \neq \gcd(X^\mathbf{u}, X^\mathbf{v}, X^\mathbf{w}),$ for every $2-$dimensional face $\{X^\mathbf{u}, X^\mathbf{v},$ $X^\mathbf{w}\} \in \nabla_\mathbf{a}.$
\item[(b)] $\gcd(X^\mathbf{u}, X^\mathbf{v})$ has $\mathcal{A}-$degree $\mathbf{a} - \mathbf{b}.$
\end{itemize}
In this case, $$\gcd(X^\mathbf{u}, X^\mathbf{v})^{-1}\Big(X^\mathbf{u} - X^\mathbf{v} \Big)$$ is the corresponding indispensable binomial of $I_\mathcal{A}.$
\end{theorem}

\begin{proof}
%
If $\mathbf{b} \prec_\mathcal{A} \mathbf{a}$ is an indispensable $\mathcal{A}-$degree of $I_\mathcal{A},$ then, by Corollary \ref{Cor Indis}, $\nabla_{\mathbf{b}} = \Big\{ \{X^{\overline{\mathbf{u}}}\}, \{X^{\overline{\mathbf{v}}}\} \Big\}.$ Let $X^\mathbf{z}$ be a monomial of $\mathcal{A}-$degree $\mathbf{a}-\mathbf{b} \neq 0$ and consider $X^\mathbf{u} = X^{\overline{\mathbf{u}}+\mathbf{z}}$ and $X^\mathbf{v} = X^{\overline{\mathbf{v}}+\mathbf{z}}.$ Notice that $\{X^\mathbf{u}, X^\mathbf{v}\} \in \nabla_\mathbf{a}$ because $\gcd(X^\mathbf{u}, X^\mathbf{v}) = X^\mathbf{z} \neq 1.$ If there exists $X^\mathbf{w} \in M_\mathbf{a} \setminus \{X^\mathbf{u}, X^\mathbf{v}\}$ with $\gcd(X^\mathbf{u}, X^\mathbf{v}, X^\mathbf{w}) = \gcd(X^\mathbf{u}, X^\mathbf{v}),$ then $X^\mathbf{z}$ divides $X^\mathbf{w}$ and so $X^\mathbf{w}/X^\mathbf{z} \in \nabla_{\mathbf{b}}.$

Conversely, let $\{X^\mathbf{u}, X^\mathbf{v}\} \in \nabla_\mathbf{a}$ be such that $\gcd(X^\mathbf{u}, X^\mathbf{v})$ has $\mathcal{A}-$degree $\mathbf{a} - \mathbf{b}$ and $\gcd(X^\mathbf{u}, X^\mathbf{v}) \neq \gcd(X^\mathbf{u}, X^\mathbf{v},$ $X^\mathbf{w}),$ for every $X^\mathbf{w} \in M_\mathbf{a} \setminus \{X^\mathbf{u}, X^\mathbf{v}\}.$ Since $\gcd(X^\mathbf{u}, X^\mathbf{v}) \neq 1,$ the monomials $X^{\mathbf{u}}/\gcd(X^\mathbf{u}, X^\mathbf{v}) $ and $X^{\mathbf{v}}/\gcd(X^\mathbf{u}, X^\mathbf{v})$ have $\mathcal{A}-$degree $\mathbf{b} \prec_\mathcal{A} \mathbf{a}$ and $\mathbf{b}$ is indispensable. Otherwise, by Corollary \ref{Cor Indis}, there exists $X^{\overline{\mathbf{w}}} \in \nabla_\mathbf{b}$ which is different from the other ones, and so $\gcd(X^\mathbf{u}, X^\mathbf{v}) = \gcd(X^\mathbf{u},$ $X^\mathbf{v}, X^\mathbf{w}),$ with $X^\mathbf{w} = X^{\overline{\mathbf{w}}} \gcd(X^\mathbf{u}, X^\mathbf{v}).$
\end{proof}

Since there exists $\mathbf{a} \in \mathcal{A}$ such that $\mathbf{b} \prec_\mathcal{A} \mathbf{a}$ for every minimal $\mathcal{A}-$degree, $\mathbf{b},$ of $I_\mathcal{A}$ (see, e.g. \cite{Sturmfels95}, Chapter 4), in order to check the existence of indispensable binomials, it is suffices to compute one (huge) simplicial complex $\nabla_\mathbf{a}$ and then use Theorem \ref{Th Indis}. Of course, this theoretical assertion is not very practical because these bounds are very coarse. Nevertheless, in some particular cases, one can find tight bounds for the minimal $\mathcal{A}-$degrees of $I_\mathcal{A}$ which combined with the high intrinsic symmetry of the simplicial complexes $\nabla_\mathbf{a}$ allows to check the existence of indispensable $\mathcal{A}-$degrees and compute all the indispensable binomials. This combined approach is applied in Section \ref{Sect An example}.

\medskip
Similar strategies may be used to compute (all) the indispensable monomials of $I_\mathcal{A}.$

\begin{definition}\label{Def IndisMon}
We say that $X^\mathbf{u} \in R_\mathbf{a}$ is an indispensable monomial of $I_\mathcal{A}$ if every system of binomial generators of $I_\mathcal{A}$ contains, at least, a binomial (up to a scalar multiple) of the form $X^\mathbf{u} - X^\mathbf{v}.$ In this case, we say that $\mathbf{a}$ is quasi-indispensable $\mathcal{A}$-degree of $I_\mathcal{A}.$
\end{definition}

Similarly to Corollary \ref{Cor Indis} we may state the following:

\begin{corollary}
Let $\mathbf{a} \in \mathcal{A}.\ \mathbf{a}$ is a quasi-indispensable if, and only if, $\nabla_\mathbf{a}$ has, at least, a $0-$dimensional connected component and $\pi^{-1}(\mathbf{a})$ has cardinality greater than or equal to $2.$
\end{corollary}

Notice that the above corollary is nothing but a combinatorial version of Theorem 3.1 in \cite{Aoki}.

It is clear that if $X^\mathbf{u} - X^\mathbf{v}$ is an indispensable binomial of $I_\mathcal{A},$ then $X^\mathbf{u}$  and $X^\mathbf{v}$ are indispensable monomials of $I_\mathcal{A}.$ Unfortunately, the converse is not true. Nevertheless, in contrast of indispensable binomials, indispensable monomials always exist (see, e.g. Proposition 3.1 in \cite{Charalambous07}).

By weakening the hypothesis in Theorem \ref{Th Indis}, we obtain the following result which allows to compute all the indispensable monomials when a sufficiently large $\mathbf{a} \in \mathcal{A}$ is known. We omit its proof because it is quite similar to the proof of Theorem \ref{Th Indis}.

\begin{corollary}
Let $\mathbf{a} \in \mathcal{A}.$ For each $X^\mathbf{u} \in M_\mathbf{a}$ and for each maximal element $\gcd(X^\mathbf{u},X^\mathbf{v})$ with respect to division in the set $\{\gcd(X^\mathbf{u},$ $X^\mathbf{w}) \mid X^\mathbf{w} \in M_\mathbf{a}\},$ the monomial $\gcd(X^\mathbf{u},X^\mathbf{v})^{-1} X^\mathbf{u}$ is an indispensable monomial of $I_\mathcal{A}.$
\end{corollary}

Notice that, in difference to Proposition 3.1 in \cite{Charalambous07}, our result does not require the previous computation of a system of generators of $I_\mathcal{A}$ to compute all its indispensable monomials.

\section{Indispensability and Gr\"obner bases}\label{Sect 03}

Let $\mathcal{A} = \mathbf{a}_1 \mathbb{N} + \ldots + \mathbf{a}_r \mathbb{N}$ be a semigroup such that $G(\mathcal{A})$ is torsion free.

\medskip
In \cite{Ohsugi} it is shown that a binomial in $I_\mathcal{A}$ is indispensable if, and only if, it or its negative belongs to the reduced Gr\"obner basis of $I_\mathcal{A}$ for any lexicographic term order on $R.$ In this section, we will prove that it is enough to check this for, at most, $r$ Gr\"obner basis with respect a degree reverse lexicographical term order on $R.$

\medskip
Fix positive integers $d_1, \ldots, d_r$ such that $I_\mathcal{A}$ is homogeneous with respect to
the grading $\deg(X_i) = d_i$ and this grading is compatible with the $\mathcal{A}-$grading of $I_\mathcal{A},$ that is to say, if $X^\mathbf{u}$ and $X^\mathbf{v}$ have the same $\mathcal{A}-$degree, then $\sum_{i=1}^r d_i u_i = \sum_{i=1}^r d_i v_i.$ This is always possible because $\mathcal{A}$ is combinatorially finite and $G(\mathcal{A})$ is torsion free.

\begin{definition}
A degree reverse lexicographic term order $\prec$ relative to the above grading on $R$ which has $X_i$ as lowest variable is any term order on $R$ represented by an $r \times r-$matrix whose first row is $(d_1, \ldots, d_r)$ and second row is $-\mathbf{e}'_i,$ where $\mathbf{e}_i$ is the $i-$th canonical basis vector of $\mathbb{Z}^r.$
\end{definition}

For a better understanding of the proof of the next result, we recall that the support of a monomial $X^\mathbf{u}$ in $R$ is $\mathrm{supp}(X^\mathbf{u}) = \mathrm{supp}(\mathbf{u}) = \big\{i \in \{1, \ldots, r\} \mid u_i \neq 0\}.$

\begin{theorem}\label{Th indisGB}
For each $i \in \{1, \ldots, r\},$ let $\prec_i$ be a degree reverse lexicographical term order on $R$ which has $X_i$ as lowest variable. A binomial in $I_\mathcal{A}$ is indispensable if, and only if, either it or its negative belongs to the reduced Gr\"obner basis $\mathcal{G}_i$ of $I_\mathcal{A}$ with respect to $\prec_i,\ i \in \{1, \ldots, r\}.$
\end{theorem}

\begin{proof}
Suppose that $X^\mathbf{u} - X^\mathbf{v} \in I_\mathcal{A}$ is indispensable. Since any reduced Gr\"obner basis of $I_\mathcal{A}$ consists of binomials (see, e.g. Proposition 1.1 in \cite{Eisenbud96}) and $X^\mathbf{u} - X^\mathbf{v}$ appears in every system of binomial generators of $I_\mathcal{A}.$ We may assume that $X^\mathbf{u} - X^\mathbf{v} \in \mathcal{G}_i,\ i = 1, \ldots r.$

Conversely, suppose that $X^\mathbf{u} - X^\mathbf{v}$ or $X^\mathbf{v} - X^\mathbf{u}$ belongs to $\mathcal{G}_i,$ for each $i \in \{1, \ldots, r\}.$ First of all, we observe that $\gcd(X^\mathbf{u}, X^\mathbf{v}) = 1.$ Otherwise, $\gcd (X^\mathbf{u},$ $X^\mathbf{v} )^{-1} \big(X^\mathbf{u}-X^\mathbf{v}\big) \in I_\mathcal{A}$ and $\gcd (X^\mathbf{u}, X^\mathbf{v} )^{-1} X^\mathbf{u}$ and $\gcd (X^\mathbf{u}, X^\mathbf{v} )^{-1} X^\mathbf{v}$ properly divide $X^\mathbf{u}$ and $X^\mathbf{v},$ respectively. Therefore, neither $X^\mathbf{u} - X^\mathbf{v}$ nor $X^\mathbf{v} - X^\mathbf{u}$ belongs to any $\mathcal{G}_i$ which is impossible by hypothesis.

Now, suppose that there exists a monomial $X^\mathbf{w}$ with the same $\mathcal{A}-$degree as $X^\mathbf{u} - X^\mathbf{v}.$ Let $X^{\overline{\mathbf{u}}} = \gcd(X^\mathbf{u}, X^\mathbf{w})$ and $X^{\overline{\mathbf{v}}} = \gcd(X^\mathbf{v}, X^\mathbf{w}).$ So, $X^\mathbf{w} = X^{\overline{\mathbf{u}}} X^{\overline{\mathbf{v}}} X^{\overline{\mathbf{w}}}$ with $\mathrm{supp}(X^{\overline{\mathbf{u}}}) \cap \mathrm{supp}(X^{\overline{\mathbf{v}}}) = \varnothing,$ because $\gcd(X^\mathbf{u}, X^\mathbf{v}) = 1.$ If $X^{\overline{\mathbf{w}}} \neq 1,$ we consider $j \in \mathrm{supp}(X^{\overline{\mathbf{w}}}).$ Then, $X^\mathbf{w} \prec_j X^\mathbf{u}$ and $X^\mathbf{w} \prec_j X^\mathbf{v},$ that is to say, $X^\mathbf{u} = \mathrm{in}_{\prec_j}(X^\mathbf{u} - X^\mathbf{w})$ and $X^\mathbf{v} = \mathrm{in}_{\prec_j}(X^\mathbf{v} - X^\mathbf{w}).$ So, by definition of Gr\"obner basis, there are two polynomials in $\mathcal{G}_j$ whose initial monomials divide $X^\mathbf{u}$ and $X^\mathbf{v},$ respectively. Therefore, neither $X^\mathbf{u} - X^\mathbf{v}$ nor $X^\mathbf{v} - X^\mathbf{u}$ could appear in $\mathcal{G}_j$ because of its reducibility. So, we may assume that $X^{\overline{\mathbf{w}}} = 1,$ that is to say, $X^\mathbf{w} = X^{\overline{\mathbf{u}}} X^{\overline{\mathbf{v}}}.$ Consider $j \in \mathrm{supp}(X^{\overline{\mathbf{v}}}).$ We have that $X^\mathbf{v} \prec_j X^\mathbf{u}.$ If $X^\mathbf{w} \prec_j X^\mathbf{v},$ then $X^\mathbf{v}$ is divisible by the initial monomial with respect to $\prec_j$ of some binomial of $I_\mathcal{A},$ thus $X^\mathbf{u} - X^\mathbf{v}$ does not belong $\mathcal{G}_j.$ Then $X^\mathbf{v} \prec_j X^\mathbf{w}$ and $X^\mathbf{v}/X^{\overline{\mathbf{v}}} \prec_j X^\mathbf{w}/X^{\overline{\mathbf{v}}} = X^{\overline{\mathbf{u}}}.$ In this case, there exists a binomial in $I_\mathcal{A}$ whose initial monomial with respect to $\prec_j$ properly divides $X^\mathbf{u} = \mathrm{in}_{\prec_j}(X^\mathbf{u} - X^\mathbf{v}).$ So we are in contradiction again.

In conclusion, $\nabla_\mathbf{a} = \big\{\{X^\mathbf{u}\},\{X^\mathbf{v}\} \big\}.$ Then, by Corollary \ref{Cor Indis}, we may assume that $X^\mathbf{u} - X^\mathbf{v}$ is an indispensable binomial of $I_\mathcal{A}.$
\end{proof}

Notice that the above theorem gives an algorithm for computing the indispensable binomials in $I_\mathcal{A}$ consisting in the computation of $r$ reduced Gr\"obner basis. Other algorithms can be found or deduced from the results in \cite{Aoki} and \cite{Charalambous07}. In the first case, $r!$ reduced Gr\"obner basis are needed. In the second case, the indispensable binomials in $I_\mathcal{A}$ are determined from one Gr\"obner basis provided that the set of minimal elements with respect $\prec_\mathcal{A}$ in the set of minimal $\mathcal{A}-$degrees of $I_\mathcal{A}$ is known.

\medskip
As immediate consequences of the above theorem we have the followings.

\begin{corollary}\label{Cor indisGB1}
For each $i \in \{1, \ldots, r\},$ let $\prec_i$ denote a degree reverse lexicographical term order on $R$ which has $X_i$ as lowest variable. If a system of generators of $I_\mathcal{A}$ is a reduced Gr\"obner basis with respect to $\prec_i,$ for every $i \in \{1, \ldots, r\},$ then $I_\mathcal{A}$ is generated by its indispensable binomials.
\end{corollary}

Observe that from the above corollary and Lemma 8.4 in \cite{Peeva2}, it follows a new proof of the result by I. Peeva and B. Sturmfels in \cite{Peeva} which states that \emph{every generic lattice ideal have a unique minimal set of binomial generators.}

\medskip
Let $\mathbf{u}_+$ and $\mathbf{u}_-$ denote the positive and negative part of $\mathbf{u} \in \mathbb{Z}^r,$ respectively. Given a system of generators $\mathcal{B}$ of $\ker(\mathcal{A}) := \ker\big(\mathbb{Z}^r \to G(\mathcal{A}); \mathbf{e}_i \mapsto \mathbf{a}_i \big),$ we write $I_\mathcal{B}$ for the binomial ideal generated by $\mathcal{G} := \big\{X^{\mathbf{u}_+} - X^{\mathbf{u}_-} \mid \mathbf{u} \in \mathcal{B} \big\}.$ Recall that $I_\mathcal{B} \subset I_\mathcal{A}$ but $I_\mathcal{B} \neq I_\mathcal{A}.$

\begin{corollary}\label{Cor indisGB2}
For each $i \in \{1, \ldots, r\},$ let $\prec_i$ denote a degree reverse lexicographical term order on $R$ which has $X_i$ as lowest variable. With the above notation, if $\mathcal{G}$ is a reduced Gr\"obner basis with respect to $\prec_i,$ for every $i \in \{1, \ldots, r\},$ then
\begin{itemize}
\item[(a)] $I_\mathcal{B} = I_\mathcal{A}$
\item[(b)] The set of indispensable binomials of $I_\mathcal{A}$ is $\mathcal{G}.$
\end{itemize}
In conclusion, $I_\mathcal{A}$ has a unique minimal system of binomial generators.
\end{corollary}

\begin{proof}
(a) By Lemma 12.2 in \cite{Sturmfels95}, $(I_\mathcal{B} : (X_1 \cdots X_r)^\infty) = I_\mathcal{A}.$ Since $\mathcal{G}$ is a reduced Gr\"obner basis with respect to a degree reverse lexicographical term order on $R$ which has $X_i$ as lowest variable, for each $i \in \{1, \ldots, r\},$ by Theorem 3.1 in \cite{Bigatti}, it follows that $(I_\mathcal{B} : (X_1 \cdots X_r)^\infty) = I_\mathcal{B}.$ So, we conclude that $I_\mathcal{B} = I_\mathcal{A}.$

Now, part (b) is an immediate consequence of the Theorem \ref{Th indisGB}.
\end{proof}

\medskip
We end this section by showing that the Lawrence ideal of $\mathcal{A}$ is generated by indispensable binomials. This result was already proved by P.Pis\'on-Casares and A. Vigneron-Tenorio in a different context (see Proposition 1(a) in \cite{PisVi}).

Recall that the Lawrence ideal of $\mathcal{A}$ is the ideal of $\mathbbmss{k}[X_1, \ldots, X_r, Y_1,$ $\ldots,$ $Y_r]$ generated by $$\Big\{X^\mathbf{u} Y^\mathbf{v} - X^\mathbf{v} Y^\mathbf{u} : \sum_{i=1}^r u_i \mathbf{a}_i = \sum_{i=1}^r v_i \mathbf{a}_i \in \mathcal{A} \Big \}.$$ Analogously, the Lawrence ideal of $\mathcal{A}$ is the ideal of $\mathbbmss{k}[X_1, \ldots, X_r,$ $Y_1, \ldots,$ $Y_r]$ associated to the Lawrence lifting of $\mathcal{A},$ that is to say, the ideal associated to the the semigroup generated by $(\mathbf{a}_1, \mathbf{e}'_1), \ldots, (\mathbf{a}_r, \mathbf{e}'_r), (0, \mathbf{e}'_1), \ldots, (0, \mathbf{e}'_r)$ in $G(\mathcal{A}) \oplus \mathbb{Z}^r,$ where $\mathbf{e}_i$ denotes the $i-$th canonical basis vector of $\mathbb{Z}^r.$

\begin{corollary}\label{Cor LI}
The Lawrence ideal of $\mathcal{A}$ is generated by indispensable binomials.
\end{corollary}

\begin{proof}
By \cite{Sturmfels95}, Theorem 7.1, any minimal binomial generating set of the Lawrence ideal of $\mathcal{A}$ form reduced Gr\"obner basis, then, by Corollary \ref{Cor indisGB1}, our claim follows.
\end{proof}

Lawrence ideals play a relevant role in the theory of toric and semigroup ideals: on the one hand, Lawrence ideals are the defining ideals toric subvarieties in a product of projective lines $\mathbb{P}^1 \times \ldots \times \mathbb{P}^1.$ On the other hand, they are used to compute the Graver basis of $I_\mathcal{A}$ (see \cite{Sturmfels95}, chapter 7, for more details).

\section{An easy example}\label{Sect An example}

In this section, we will apply our results to prove that the toric ideal associated to
the binary marginal independence model $\mathcal{A}$ induced by the undirected graph $G = (V,E)$ with
$V = \big\{ \{1\}, \{2\},$ $\{3\}, \{4\} \big\}$ and $E =
\big\{ \{1, 2\}, \{2, 3\}, \{3, 4\}, \{4, 1\}, \{2, 4\} \big\},$
is generated by indispensable binomials.

We refer to the interested reader to \cite{Geiger} for the details omitted here about graphical models.

\medskip
The model $\mathcal{A}$ has associated matrix
$$
A = {\tiny \left(
\begin{array}{*{16}{r}}
1 & 0 & 0 & 0 & 0 & 0 & 0 & 0 & 1 & 0 & 0 & 0 & 0 & 0 & 0 & 0\\
0 & 1 & 0 & 0 & 0 & 0 & 0 & 0 & 0 & 1 & 0 & 0 & 0 & 0 & 0 & 0\\
0 & 0 & 1 & 0 & 0 & 0 & 0 & 0 & 0 & 0 & 1 & 0 & 0 & 0 & 0 & 0\\
0 & 0 & 0 & 1 & 0 & 0 & 0 & 0 & 0 & 0 & 0 & 1 & 0 & 0 & 0 & 0\\
0 & 0 & 0 & 0 & 1 & 0 & 0 & 0 & 0 & 0 & 0 & 0 & 1 & 0 & 0 & 0\\
0 & 0 & 0 & 0 & 0 & 1 & 0 & 0 & 0 & 0 & 0 & 0 & 0 & 1 & 0 & 0\\
0 & 0 & 0 & 0 & 0 & 0 & 1 & 0 & 0 & 0 & 0 & 0 & 0 & 0 & 1 & 0\\
0 & 0 & 0 & 0 & 0 & 0 & 0 & 1 & 0 & 0 & 0 & 0 & 0 & 0 & 0 & 1\\
1 & 0 & 1 & 0 & 0 & 0 & 0 & 0 & 0 & 0 & 0 & 0 & 0 & 0 & 0 & 0\\
0 & 1 & 0 & 1 & 0 & 0 & 0 & 0 & 0 & 0 & 0 & 0 & 0 & 0 & 0 & 0\\
0 & 0 & 0 & 0 & 1 & 0 & 1 & 0 & 0 & 0 & 0 & 0 & 0 & 0 & 0 & 0\\
0 & 0 & 0 & 0 & 0 & 1 & 0 & 1 & 0 & 0 & 0 & 0 & 0 & 0 & 0 & 0\\
0 & 0 & 0 & 0 & 0 & 0 & 0 & 0 & 1 & 0 & 1 & 0 & 0 & 0 & 0 & 0\\
0 & 0 & 0 & 0 & 0 & 0 & 0 & 0 & 0 & 1 & 0 & 1 & 0 & 0 & 0 & 0\\
0 & 0 & 0 & 0 & 0 & 0 & 0 & 0 & 0 & 0 & 0 & 0 & 1 & 0 & 1 & 0\\
0 & 0 & 0 & 0 & 0 & 0 & 0 & 0 & 0 & 0 & 0 & 0 & 0 & 1 & 0 & 1
\end{array}
\right)}
$$
that, in a more condensed, can be written as
$$A = \left(\begin{array}{c} 1'_2 \otimes I_4 \otimes I_2 \\ I_4 \otimes 1'_2 \otimes I_2 \end{array}\right) \in \mathbb{Z}^{16 \times 16},$$ where $1'_2 = (1\ 1), I_n$ is the $n \times n-$identity matrix and the symbol $\otimes$ denotes the Kronecker product.

Let $I_\mathcal{A} \subset \mathbbmss{k}[X_1, \ldots, X_{16}]$ be the semigroup ideal associated to the subsemigroup of $\mathbb{N}^{16}$ generated by the columns of $A.$ That is to say, $I_\mathcal{A}$ is the toric ideal of the independence model $\mathcal{A}.$

The reader may note that by performing row operations on $A,$ it becomes a [Lawrence lifting]-type matrix, so, by Corollary \ref{Cor LI}, we may conclude that $I_\mathcal{A}$ is generated by indispensable binomials and hence it has a unique minimal system of binomial generators. Nevertheless, in order to illustrate the results in this paper, we will proceed by assuming that we do not know this fact.

\medskip
Since $\mathcal{A}$ is a decomposable graphical model, $I_\mathcal{A}$ has a quadratic Gr\"obner basis (see Theorem 4.3 in \cite{Geiger}). Therefore, any indispensable binomial (if exist) is homogeneous (with the usual grading) and has total degree equals two. Moreover, every quadratic binomial in $I_\mathcal{A}$ consists in differences of square-free monomials, because no sum of two (not necessarily different) columns of $A$ is equal to the double of another one. In conclusion, if $X^{\mathbf{u}} - X^{\mathbf{v}}$ is an indispensable binomial in $I_\mathcal{A},$ then $u_i \leq 1$ and $v_i \leq 1,$ for all $i,$ where $\mathbf{u} = (u_1, \ldots, u_r)'$ and $\mathbf{v} = (v_1, \ldots, v_r)' \in \mathbb{N}^r,$ as usual.

Thus, the $\mathcal{A}-$degree $\mathbf{a} =  A\, \mathbf{e},$ with $\mathbf{e} = (1, \ldots, 1)' \in \mathbb{Z}^{16},$ ``captures'' all the indispensable binomials.

Now, let us compute the non-negative integer solutions of the system $A \mathbf{u} = \mathbf{a}.$ The general solution of the linear system of equation $A \mathbf{u} = \mathbf{a}$ over $\mathbb{Q}$ is $$\begin{array}{ccccc} \mathbf{u} = & \big(\ 1+a, & 1+b, & 1-a, & 1-b,\\ & \mbox{}\ 1+c, & 1+d, & 1-c, & 1-d,\\ &\mbox{}\ 1-a, & 1-b, & 1+a, & 1+b,\\ &\mbox{}\ 1-c, & 1-d, & 1+c, & 1+d\ \big)'
\end{array}$$ Thus, it is clear that the non-negative integer solutions correspond to the values of $a,b,c$ and $d$ in $\{-1,0,1\}.$

First of all we observe that any $\mathbf{u} \in \mathbb{N}^{16}$ such that $A \mathbf{u} = \mathbf{a}$ has the form $$(\mathbf{p}'_i , \mathbf{p}'_j, \mathbf{p}'_{\sigma(i)}, \mathbf{p}'_{\sigma(j)})'$$ with
\begin{equation}\label{ecu pes}
\begin{array}{c|cccc|cccc}
\mathbf{p}_1 & \mathbf{p}_2 & \mathbf{p}_3 &
\mathbf{p}_4 & \mathbf{p}_5 & \mathbf{p}_6 &
\mathbf{p}_7 & \mathbf{p}_8 & \mathbf{p}_9 \\
\cline{1-9}
1 & 0 & 1 & 2 & 1 & 0 & 2 & 2 & 0\\
1 & 1 & 0 & 1 & 2 & 0 & 0 & 2 & 2\\
1 & 2 & 1 & 0 & 1 & 2 & 0 & 0 & 2\\
1 & 1 & 2 & 1 & 0 & 2 & 2 & 0 & 0
\end{array}
\end{equation}
and $\sigma = (24)(35)(68)(79).$

Therefore, it follows that $\nabla_\mathbf{a}$ is a direct product of simplicial complexes. Concretely, $\nabla_\mathbf{a} \cong K \times K,$ with $$K = \Big\{F \subseteq \{Y^{\mathbf{p}_i} \mid i = 1, \ldots, 9\}\ \mid\ \gcd(F) \neq 1 \Big\}.$$ Furthermore, given $\mathbf{u} = (\mathbf{p}'_i , \mathbf{p}'_j, \mathbf{p}'_{\sigma(i)}, \mathbf{p}'_{\sigma(j)})'$ and $\mathbf{v} = (\mathbf{p}'_k , \mathbf{p}'_l, \mathbf{p}'_{\sigma(k)}, \mathbf{p}'_{\sigma(l)})',$ by Corollary \ref{Cor Indis} and Theorem \ref{Th Indis}, we have that $\gcd(X^\mathbf{u}, X^\mathbf{v})^{-1}\Big(X^\mathbf{u} - X^\mathbf{v} \Big)$ is indispensable if, and only if, $$i=k,\ j \neq l\ \text{and}\ \mathrm{gcd}(Y^{\mathbf{p}_j}, Y^{\mathbf{p}_l})\ \text{is uniquely attained}$$ or $$i \neq k,\ j=l\ \text{and}\ \mathrm{gcd}(Y^{\mathbf{p}_i}, Y^{\mathbf{p}_k})\ \text{is uniquely attained}.$$ In fact, the indispensable binomials are \begin{equation}\label{ecu indis}\mathrm{gcd}\big(Y^{\mathbf{p}_i}Z^{\mathbf{p_{\sigma(i)}}}, Y^{\mathbf{p}_j}Z^{\mathbf{p_{\sigma(j)}}}\big)^{-1} \Big(Y^{\mathbf{p}_i}Z^{\mathbf{p_{\sigma(i)}}} - Y^{\mathbf{p}_j}Z^{\mathbf{p_{\sigma(j)}}} \Big)\end{equation} with $\mathrm{gcd}(Y^{\mathbf{p}_i}, Y^{\mathbf{p}_j})$ uniquely attained and $Y_k = X_k, Z_k = X_{k+8},\ k \in \{1, \ldots, 4\}$ or $Y_k = X_k, Z_k = X_{k+8},\ k \in \{5, \ldots, 8\}.$

Summarizing, the indispensable binomials of $I_\mathcal{A}$ are determined by the pair of vertices of $K,\ \mathcal{P},$ whose greatest common divisor is different from the greatest common divisor of any $2-$dimensional face of $K$ containing them.

\medskip
Notice that the natural action of $H = \langle (1\ 2\ 3\ 4) \rangle$ on $\{Y_1, Y_2, Y_3, Y_4\}$ leaves $K$ invariant. In fact, this action is the same than the one given by $\tilde H = \langle (2\ 3\ 4\ 5)(6\ 7\ 8\ 9) \rangle$ on the set of $\mathbf{p}_i$'s. Thus, in order to compute $\mathcal{P},$ it suffices to perform the computation modulo $H.$ There are nine $1-$dimensional faces different modulo $H:$
$$\{Y^{\mathbf{p_1}}, Y^{\mathbf{p_2}}\}, \{Y^{\mathbf{p_1}}, Y^{\mathbf{p_6}}\}, \{Y^{\mathbf{p_2}}, Y^{\mathbf{p_3}}\}, \{Y^{\mathbf{p_2}}, Y^{\mathbf{p_4}}\}, \{Y^{\mathbf{p_2}}, Y^{\mathbf{p_6}}\}$$
$$\{Y^{\mathbf{p_2}}, Y^{\mathbf{p_7}}\}, \{Y^{\mathbf{p_2}}, Y^{\mathbf{p_8}}\}, \{Y^{\mathbf{p_2}}, Y^{\mathbf{p_9}}\}, \{Y^{\mathbf{p_6}}, Y^{\mathbf{p_7}}\},$$ but only three of them have greatest common divisor different from any $2-$dimensional face of $K:$ $$\{Y^{\mathbf{p_1}}, Y^{\mathbf{p_2}}\}, \{Y^{\mathbf{p_2}}, Y^{\mathbf{p_6}}\}, \{Y^{\mathbf{p_2}}, Y^{\mathbf{p_9}}\}.$$ An easy computation shows that
$$\mathrm{gcd}\big(Y^{\mathbf{p}_i}, Y^{\mathbf{p}_j}\big)^{-1} \Big(Y^{\mathbf{p}_i} - Y^{\mathbf{p}_j} \Big)
= \left\{\begin{array}{r} \pm(Y_1 - Y_3) \smallskip \\ \pm(Y_2 - Y_4) \end{array}\right.$$ when $\mathrm{gcd}(Y^{\mathbf{p}_i}, Y^{\mathbf{p}_j})$ is uniquely attained.

Therefore, by (\ref{ecu indis}), we conclude that there are four indispensable binomials in $I_\mathcal{A}:$
$$X_1X_{11}-X_3X_9,\ X_2X_{12}-X_4X_{10},\ X_5X_{15}-X_{7}X_{13}\quad \text{and}\quad X_6X_{16}-X_8X_{14}.$$
Moreover, since these four binomials form a Gr\"obner basis with respect to any term order on $R$ (because
their initial monomials have disjoint support) and their exponent vectors generate $\ker(\mathcal{A}),$ by Corollary \ref{Cor indisGB2}, we may assure that $I_\mathcal{A}$ is generated by its indispensable binomials, that is to say, the toric ideal associated to the binary marginal independence model $\mathcal{A}$ has a unique minimal system of binomial generators.

\medskip \noindent
\textbf{Acknowledgments.-}
We want to warmly thank the anonymous referee for his/her comments and remarks.


\begin{thebibliography}{12}

\bibitem{Aoki}
\textsc{S. Aoki, A. Takemura, R. Yoshida.}
\newblock \emph{Indispensable monomials of toric ideals and Markov bases.}
J. Symbolic Comput. \textbf{43} (2008), no. 6-7, 490--507.

\bibitem{Bigatti}
\textsc{A. Bigatti, R. {La Scala}, L. Robbiano.}
\newblock \emph{Computing toric ideals.}
\newblock J. Symbolic Comput. \textbf{27} (1999), 351--365.

\bibitem{Collectanea}
\textsc{E. Briales, A. Campillo, C. Mariju\'an, P. Pis\'on},
\newblock \emph{Combinatorics of syzygies for semigroup algebra.}
\newblock {Collect. Math.} \textbf{49} (1998), 239--256.

\bibitem{BCMP}
\textsc{E. Briales, A. Campillo, C. Mariju\'an, P. Pis\'on.}
\newblock \emph{Minimal Systems of Generetors for Ideals of Semigroups.}
\newblock  J. Pure Appl. Algebra, \textbf{124} (1998), 7--30.

\bibitem{Charalambous07}
\textsc{H. Charalambous, A. Katsabekis, A. Thoma.}
\newblock \emph{Minimal systems of binomial generators and the indispensable complex of a toric ideal.}
\newblock Proc. Amer. Math. Soc. \textbf{135} (2007), 3443--3451.

\bibitem{Charalambous08}
\textsc{H. Charalambous, A. Thoma.}
\newblock \emph{On simple $\mathcal{A}-$multigraded minimal resolutions.}
\newblock Preprint (2008).

\bibitem{Diaconis}
\textsc{P. Diaconis, B. Sturmfels.}
\newblock \emph{Algebraic algorithms for sampling from conditional distributions.}
\newblock Ann. Statist. \textbf{26}(1) (1998), 363--397.

\bibitem{Eisenbud96}
\textsc{D. Eisenbud, B. Sturmfels}
\newblock \emph{Binomial ideals.}
\newblock Duke Math. J. 84 (1996), no. 1, 1--45.

\bibitem{Eliahou}
\textsc{S. Eliahou.}
\newblock \emph{Courbes monomiales et alg\'{e}bre de Rees symbolique.}
\newblock PhD Thesis. Universit\'{e} of Gen\`{e}ve, 1983.

\bibitem{Geiger}
\textsc{D. Geiger, C. Meek, B. Sturmfels.}
\newblock \emph{On the toric algebra of graphical models.}
\newblock Ann. Statist. 34 (2006), no. 3, 1463--1492.

\bibitem{Herzog70}
\textsc{J. Herzog.}
\newblock \emph{Generators and relations of abelian semigroups and semigroup rings.}
\newblock Manuscripta Math. 3 1970 175--193.

\bibitem{Miller05}
\textsc{E. Miller, B. Sturmfels}
\newblock \emph{Combinatorial Commutative Algebra}.
\newblock Vol. 227 of Graduate Texts in Mathematics. Springer, New York. 2005.

\bibitem{Ohsugi05}
\textsc{H. Ohsugi, T. Hibi}
\newblock \emph{Indispensable binomials of finite graphs.}
\newblock J. Algebra Appl. 4 (2005), no. 4, 421--434.

\bibitem{Ohsugi}
\textsc{H. Ohsugi, T. Hibi.}
\newblock \emph{Toric ideals arising from contingency tables,}
\newblock in Commutative Algebra and Combinatorics.
\newblock Ramanujan Mathematical Society Lecture Notes Series, Vol. 4, Ramanujan Mathematical Society, Mysore, India, 2007, pp. 91-115.

\bibitem{OjVi}
\textsc{I. Ojeda, A. Vigneron-Tenorio.}
\newblock \emph{Simplicial complexes and minimal free resolution of monomial algebras.}
\newblock \url{http://dx.doi.org/10.1016/j.jpaa.2009.08.009}

\bibitem{OjPis}
\textsc{I. Ojeda , P. {Pis\'on-Casares}.}
\newblock \emph{On the hull resolution of an affine monomial curve.}
\newblock J. Pure Appl. Algebra \textbf{192} (2004), 53--67.

\bibitem{Oj2}
\textsc{I. Ojeda.}
\newblock \emph{Examples of generic lattice ideals of codimension 3.}
\newblock Comm. Algebra 36 (2008) 279-287.

\bibitem{Peeva}
\textsc{I. Peeva, B. Sturmfels.}
\newblock \emph{Generic lattice ideals},
\newblock J. Amer. Math. Soc. \textbf{11} (1998), 363--373.

\bibitem{Peeva2}
\textsc{I. Peeva, B. Sturmfels.}
\newblock \emph{Syzygies of codimension 2 lattice ideals},
\newblock Math. Z. \textbf{229}, (1998), 163–-194.

\bibitem{PisVi}
\textsc{P. Pis\'on-Casares, A. Vigneron-Tenorio.}
\newblock \emph{On Lawrence semigroups},
\newblock J. Symbolic Comput. \textbf{43} (2008), 804--810.


\bibitem{Rosales}
\textsc{J.C. Rosales, P.A. Garc\'{\i}a-S\'anchez}
\newblock \emph{Finitely generated commutative monoids},
\newblock Nova Science Publishers, Inc., New York, 1999.

\bibitem{Sturmfels95}
\textsc{B. Sturmfels.}
\newblock \emph{Gr\"obner bases and convex polytopes},
\newblock volume~8 of \emph{University  Lecture Series.}
\newblock American Mathematical Society, Providence, RI, 1996.

\bibitem{Takemura}
\textsc{A. Takemura, S. Aoki.}
\newblock \emph{Some characterizations of minimal Markov basis for sampling from discrete conditional distributions.} \newblock Ann. Inst. Statist. Math. \textbf{56}(1) (2004), 1--17.

\end{thebibliography}
\end{document}